\documentclass[english]{amsart}

\usepackage{babel}
\usepackage{xypic}
\usepackage{amstext}
\usepackage{amsmath}
\usepackage{amsfonts}
\usepackage{latexsym}
\xyoption{all}

\begin{document}

\newcommand{\id}{{\rm id}}
\newcommand{\sO}{{\mathcal O}}
\newcommand{\I}{{\mathcal I}}
\newcommand{\J}{{\mathcal J}}
\newcommand{\G}{{\mathbb G}}
\newcommand{\rk}{{\rm rk}}
\newcommand{\codim}{{\rm codim}}
\newcommand{\Pic}{{\rm Pic}}
\newcommand{\GHom}{{\mathcal H}om}
\newcommand{\GExt}{{\mathcal E}xt}
\newcommand{\Hom}{{\rm Hom}}
\newcommand{\Ext}{{\rm Ext}}
\newcommand{\End}{{\mathcal E}nd}
\newcommand{\PN}{{\mathbb P}}
\newcommand{\lra}{\longrightarrow}
\newcommand{\KC}{{\mathbb C}}
\newcommand{\ind}{{\rm ind}}
\newcommand{\im}{{\rm im}}
\newcommand{\Seq}{{\mathcal S}}

\renewcommand{\theequation}{\arabic{section}.\arabic{lemma}.\arabic{equation}}
\newcounter{lemma}
\renewcommand{\thelemma}{\arabic{section}.\arabic{lemma}}
\newtheorem{lemma1}[lemma]{\setcounter{equation}{0}}

\newenvironment{lemma}{\begin{lemma1}{\bf Lemma.}}{\end{lemma1}}
\newenvironment{example}{\begin{lemma1}{\bf Example.}\rm}{\end{lemma1}}
\newenvironment{abs}{\begin{lemma1}\rm}{\end{lemma1}}
\newenvironment{theorem}{\begin{lemma1}{\bf Theorem.}}{\end{lemma1}}
\newenvironment{theorem2}[1]{\begin{lemma1}{\bf Theorem [#1].}}{\end{lemma1}}
\newenvironment{proposition}{\begin{lemma1}{\bf Proposition.}}{\end{lemma1}}
\newenvironment{corollary}{\begin{lemma1}{\bf Corollary.}}{\end{lemma1}}
\newenvironment{remark}{\begin{lemma1}{\bf Remark.}\rm}{\end{lemma1}}
\newenvironment{definition}{\begin{lemma1}{\bf Definition.}}{\end{lemma1}}

\renewcommand{\labelenumii}{(\roman{enumii})}


\title{Submanifolds with splitting tangent sequence}
\author{Priska Jahnke}
\address{Mathematisches Institut \\ Universit\"at Bayreuth \\ D--95440 Bayreuth/Germany}
\email{priska.jahnke@uni-bayreuth.de}
\date{\today}
\thanks{The author was supported by the Deutsche Forschungsgemeinschaft.}
\maketitle

\section*{Introduction}

In the late 50's A. Van de Ven stated the following theorem (\cite{VdV}): the only compact submanifolds with splitting tangent sequence of the projective space are linear subspaces. For a submanifold $X$ of a projective manifold $M$ the {\em tangent sequence} denotes the natural sequence  
 \[0 \lra T_X \lra T_M|_X \lra N_{X/M} \lra 0\]
mapping the tangent bundle of $X$ into the restriction of $T_M$ to $X$. If the sequence splits, $T_X$ becomes a direct summand of $T_M|_X$. This means that the first linear approximations of $X$ and $M$ coincide in a neighborhood of $X$. The linear structure of $\PN_n$ then forces $X$ to be a linear subspace. Most proofs of Van de Ven's theorem rely on this fact.

\vspace{0.2cm}

Aim of this paper is to study submanifolds with splitting tangent sequence of a given projective manifold $M$. We will consider in particular homogeneous manifolds, for example quadrics or abelian varieties. A compact manifold $M$ is homogeneous if and only if its tangent bundle is globally generated. Therefore any submanifold with splitting tangent sequence of a homogeneous manifold is again homogeneous.

The idea is to study positivity of the tangent bundle: if $T_M$ is ``positive'', then the quotient $T_X$ should be ``positive'' as well, telling us in turn hopefully something about the geometry of $X$. The main reason why this works out for homogeneous manifolds is the theorem of Borel and Remmert (\cite{BoRe}): if $M$ is a homogeneous manifold, then $M$ decomposes as
 \[M \simeq A \times G/P,\] 
where $A$ is an abelian variety (corresponding to the ``trivial part'' of the tangent bundle) and $G/P$ is rational homogeneous (the ``positive'' part of $T_M$). Rational homogeneous manifolds are quotients of semi-simple complex Lie groups and therefore completely classified in terms of Dynkin diagrams. For the abelian part we can show

\vspace{0.2cm}

\noindent {\bf \ref{Tori}~Theorem. }{\em The tangent sequence of a submanifold of an abelian variety splits if and only if it is abelian.}

\vspace{0.2cm}

The question arises, whether the analogous statement holds true for the positive part as well: {\em Is a submanifold with splitting tangent sequence of a rational homogeneous manifold again rational homogeneous?} So far we cannot deal with this question in general, but the answer turns out to be ``yes'' for all $M$ considered here (and it is ``yes'' provided that the dimension of $X$, relative to some positivity of $T_M$ is big enough, see Proposition~\ref{pos}).

\vspace{0.2cm}

\noindent {\bf \ref{Qn}~Theorem. }{\em The tangent sequence of a submanifold of dimension at least $2$ of the quadric $Q_n$ splits, if and only if it is a linear subspace or a complete intersection subquadric. A curve in $Q_n$ with splitting tangent sequence is rational.}

\vspace{0.1cm}

\noindent Key of the proof is the fact that $Q_n$ and the projective space are the only projective manifolds $X$ such that $\bigwedge^2T_X$ is ample. This is a result of Cho and Sato (\cite{ChoSato}). But in our case we know moreover, that $X$ is homogeneous, meaning that in fact we do not need that strong result, just using classification (which is of course in no way easier). Next consider Grassmannians: 

\vspace{0.2cm} 

\noindent {\bf \ref{homratG}~Proposition. }{\em If $X$ is a submanifold with splitting tangent sequence of a Grassmannian, then $X$ is rational homogeneous.}

\vspace{0.2cm}

One more homogeneous example is the flag manifold $\PN(T_{\PN_2})$. In Proposition~\ref{PvT} we show, that the only submanifolds with splitting tangent sequence of $\PN(T_{\PN_2})$ are sections of lines in $\PN_2$ corresponding to one of the two projections $\pi_i: \PN(T_{\PN_2}) \to \PN_2$. Finally we consider the blow-up $\PN_n(p)$ of $\PN_n$ in a point $p$, an almost homogeneous manifold, obtaining a very similar result:

\vspace{0.2cm}

\noindent {\bf \ref{Pn(p)}~Theorem. }{\em The tangent sequence of a submanifold $X$ in $\PN_n(p)$ splits, if and only if one of the following holds: (i) $X$ is a linear subspace of the exceptional divisor, (ii) $X$ is the strict transform of a linear subspace in $\PN_n$, (iii) $X$ is the strict transform of a smooth conic through $p$ or of a plane curve of degree $3$ with node in $p$.}

\vspace{0.2cm}

\noindent {\bf Acknowledgements.} This article contains results from the author's thesis which was prepared at the graduate program ``Complex Manifolds'' at the University of Bayreuth. The author wants to express her gratitude to her adviser Th. Peternell and to the late M. Schneider. As well she would like to thank the referee for his very careful proof-reading and his many valuable remarks.


\section{Splitting of sequences}
\setcounter{lemma}{0}

After fixing some notations we will sum up some criterions for submanifolds with splitting tangent sequence. We will need these properties frequently in the sequel, so it might be convenient for the reader to collect them here.

\begin{abs}
{\bf Notations and Conventions.} All manifolds (varieties) are assumed to be projective and defined over $\KC$. Let $X \subset M$ be a submanifold. Then the {\em tangent sequence} of $X$ in $M$
\newcommand{\myequation}{\arabic{section}.\arabic{equation}} \renewcommand{\theequation}{$\Seq_{X/M}$}\begin{equation}
   0 \lra T_X \lra T_M|_X \lra N_{X/M} \lra 0
 \end{equation}\renewcommand{\theequation}{\myequation}\addtocounter{equation}{-1}
will be denoted by $\Seq_{X/M}$. If $Z \subset X$ is a further submanifold, then we have another exact sequence
  \[0 \lra N_{Z/X} \lra N_{Z/M} \lra N_{X/M}|_Z \lra 0,\]
the {\em normal bundle sequence} of $Z,X,M$. By definition, a short exact sequence of vector bundles $0 \to A \to B \stackrel{\alpha}{\to} C \to 0$ splits, if there exists a morphism $\beta:C \to B$ with $\alpha \circ\beta = \id_C$. This is true if $\Ext^1(C, A) = 0$, called {\rm $\Ext$}- or {\em $h^1$-criterion}. The Picard number of a manifold $X$ will be denoted by $\rho(X)$, the canonical divisor by $K_X$. We will identify line bundles and divisors and write $K_X$ instead of $\sO_X(K_X)$. If $X$ is rational homogeneous, the (positive) generators of $\Pic(X)$ are called {\em fundamental divisors}; if $\rho(X) = 1$, the fundamental divisor is denoted by $\sO_X(1)$. If $E$ is a vector bundle on $X$, then $\PN(E)$ is the projective bundle of hyperplanes in the sense of Grothendieck, ${\bf P}(E)$ is the projective bundle of lines.
\end{abs}

\begin{abs} \label{trans}
{\bf Transitivity.} Let $Z \subset X$ and $X \subset M$ be submanifolds. Then we have
\[\xymatrix{&&& 0 \ar[d]&\\
  0 \ar[r] & T_Z \ar[r] \ar@{=}[d] & T_X|_Z \ar[r] & N_{Z/X} \ar[r] \ar[d] & 0\\
  0 \ar[r] & T_Z \ar[r] & T_M|_Z \ar[r] \ar@{=}[d] & N_{Z/M} \ar[r] \ar[d] & 0\\
  0 \ar[r] & T_X|_Z \ar[r]  & T_M|_Z \ar[r]  & N_{X/M}|_Z \ar[r] \ar[d] & 0\\
  &&& 0 &}\]
We obtain by chasing the splitting obstruction in that diagram (or by constructing the splitting morphism directly):
\begin{enumerate}
 \item If $\Seq_{Z/M}$ splits, that is the (horizontal) line in the middle, then $\Seq_{Z/X}$ splits, that is the line on top.
 \item If $\Seq_{Z/X}$ and $\Seq_{X/M}$ both split, that are the top and the bottom (horizontal) lines, then the (vertical) normal bundle sequence and $\Seq_{Z/M}$ split. Note that it is enough that $\Seq_{X/M}$ splits on $Z$.
\end{enumerate}
\end{abs}

\begin{abs} \label{fiberbundle}
{\bf Fiber bundles.} Let $M$ be a projective manifold and $\pi:M \to Y$ a generic fiber bundle with typical fiber $X$ onto some normal reduced variety $Y$. Then $\Seq_{X/M}$ splits. To see this consider an (analytically) open set $U \subset Y$, where the fiber bundle is trivial and $Y$ is smooth. Then the normal bundle of $X$ in the preimage $\pi^{-1}(U)$ is trivial. On the other hand, the tangent bundle of $\pi^{-1}(U)$ is the restriction of the tangent bundle of $M$. We may hence assume that $Y$ is smooth and $M = Y \times X$ is a product. In this case, the assertion is trivial.
\end{abs}

\begin{abs} \label{projbundle}
{\bf Projective bundles.} Let $E$ be a vector bundle on the projective manifold $M$ and consider the projective bundle $\pi: \PN(E) \to M$. Let $X \subset \PN(E)$ be a submanifold with splitting tangent sequence, such that the image $\pi(X)$ is smooth. Assume 
 \begin{enumerate} 
  \item $X = \PN(E|_{\pi(X)})$ or
  \item $\pi|_X:X \to \pi(X)$ is an isomorphism and $\Hom(T_{\PN(E)/M}|_X, T_X) = 0$.
 \end{enumerate}
Then $\Seq_{\pi(X)/M}$ splits. To prove this consider the diagram
\[\xymatrix{&& 0 \ar[d] & 0 \ar[d]&\\
  0 \ar[r] & T_{X/\pi(X)} \ar[r]^{\alpha_1} \ar[d]_{\gamma_1} & T_X \ar[r]^{\alpha_2} \ar[d]_{\gamma_2} & \pi^*T_{\pi(X)} \ar[r] \ar[d]_{\gamma_3} & 0\\
  0 \ar[r] & T_{\PN(E)/M}|_X \ar[r]^{\beta_1} & T_{\PN(E)}|_X \ar[r]^{\beta_2} \ar[d] \ar@/_/[u]_{\psi} & \pi^*T_M|_{\pi(X)} \ar[r] \ar[d] \ar@{..>}@/_/[u]_{\varphi} & 0\\
  & & N_{X/\PN(E)} \ar[r] \ar[d] & \pi^*N_{\pi(X)/M} \ar[d] & \\
  && 0 & 0 &}\]
In the first case we have $\gamma_1 = \id$; in the second case the relative tangent bundle $T_{X/\pi(X)}$ vanishes, hence $\gamma_1 = 0$. The additional assumption in this case implies $\psi \circ \beta_1 = 0$. In both cases we construct the splitting map $\varphi$ as follows: for an element $x \in \pi^*T_M|_{\pi(X)}$ there exists some $y \in T_{\PN(E)}|_X$ with $\beta_2(y) = x$. Define 
 \[\varphi(x) := \alpha_2(\psi(y)).\] 
Since $y$ is not unique, we of course have to check that $\varphi$ is well-defined. This, and the splitting property $\varphi \circ \gamma_3 = \id_{\pi^*T_{\pi(X)}}$ can be shown by diagram chase. Since $\pi|_X$ has connected fibers, we conclude the splitting of $\Seq_{\pi(X)/M}$ from the splitting of the right (vertical) column.

\vspace{0.2cm}

Assume conversely that $X \subset \PN(E)$ is a submanifold, such that $\pi|_X: X \to \pi(X)$ is an isomorphism and $\Seq_{\pi(X)/M}$ splits. Then $\Seq_{X/\PN(E)}$ splits. Indeed: the splitting of $\Seq_{\pi(X)/M}$ implies the existence of a splitting map $\varphi$ in the above diagram, and we have to show the existence of $\psi$. Since $X$ is a section, $T_{X/\pi(X)}$ is zero, hence $\alpha_2$ is an isomorphism. We define $\psi = \alpha_2^{-1} \circ \varphi \circ \beta_2$. 
\end{abs}


\section{Positivity of the tangent bundle} \label{sec pos}
\setcounter{lemma}{0}

Throughout this section let $M = G/P$ be a rational homogeneous manifold of dimension $n$. Then $M$ is in particular a Fano manifold; we denote its index by $r$, i.e. 
 \[-K_M = rH,\] 
where $H \in \Pic(M)$ is ample and $r$ is chosen maximal with this property. The tangent bundle of any homogeneous manifold $M$ is globally generated; we denote the associated morphism by $\Phi_M$,
 \begin{equation} \label{amp}
   \Phi_M: \PN(T_M) \lra \PN(H^0(M, T_M)).
 \end{equation}
After Stein factorization we may assume that the fibers of $\Phi_M$ are connected. If the maximal dimension of a fiber of $\Phi_M$ is $k$, then $T_M$ is called $k$-ample by Sommese (\cite{Sommese}). This is equivalent to saying that $k$ is the maximal dimension of a subvariety of $M$, where the restriction of $T_M$ has a trivial quotient. Note that $k$-ample vector bundles behave just like ample ones (see \cite{Sommese}): quotients of a $k$-ample vector bundle are again $k$-ample, and there is a cohomological criterion for $k$-ample bundles as well. Goldstein computed the co-$k$-ampleness ${\rm ca}(M) = n-k$ for any rational homogeneous manifold $M$ in \cite{Gold}. 

\vspace{0.2cm}

Another way to characterize positivity is to study exterior powers of $T_M$. For $M \not= \PN_n$, the tangent bundle itsself is not ample, but its determinant $-K_M = \bigwedge^n T_M$ is. We may hence ask for the minimal $l$, such that $\bigwedge^lT_M$ is ample. Note that this implies that $\bigwedge^jT_M$ is ample for all $j \ge l$. By Mori's theorem, $l = 1$ implies $M \simeq \PN_n$, and Cho and Sato have proved in \cite{ChoSato} that $l = 2$ holds true if and only if $M \simeq Q_n$ (for $n \ge 3$). For general (non-homogeneous) Fano manifolds it is not known, if it is possible to give an estimate for $l$ in terms of index and dimension of $M$ (see \cite{Pet}, (2.12.)). But for homogeneous $M$ we will see that it is not difficult to compute $l$. First we prove the following easy, but nevertheless useful statement:

\begin{proposition} \label{pos}
Let $M$ be rational homogeneous and $X \subset M$ a submanifold with splitting tangent sequence. Let $k, l$ be the minimal numbers, such that $T_M$ is $k$-ample and $\bigwedge^lT_M$ is ample, respectively. If $\dim(X) \ge \min(k+1, l)$, then $X$ is rational homogeneous. 
\end{proposition}

\begin{corollary}
In the situatuion of the Proposition, let $X \simeq A \times G/P$ be the Borel-Remmert decomposition of $X$. Then $\dim(A) \le \min(k, l-1)$.
\end{corollary}

\begin{proof}[Proof of Proposition~\ref{pos}]
Assume first $\dim(X) \ge k+1$. The splitting $T_M|_X \to T_X \to 0$ yields that $T_X$ is $k$-ample and $X$ is homogeneous. Let $X  \simeq A \times G/P$ be the Borel-Remmert decomposition of $X$ and let $p_A:X \to A$ be the projection. Then there is a surjective map $T_X \to p_A^*T_A$. This is not possible if $\dim(A) > 0$, since $T_A$ is trivial and $k < \dim(X)$.

If $\dim(X) \ge l$, then $\bigwedge^lT_X$ is ample, since the splitting gives surjective maps $\bigwedge^jT_M|_X \to \bigwedge^jT_X \to 0$ for all $j \le \dim(X)$. Hence $-K_X$ is ample. 
\end{proof} 

\begin{remark} \label{kl}
For rational homogeneous $M$, the numbers $k,l$ are related to the index $r$ as follows: Goldstein computed ${\rm ca}(M) = n -k$. Comparing his results with $r$, for example found in \cite{Akh}, we obtain
 \[k \le n - r +1.\]
Equality holds for example if $\rho(M) = 1$ and the Lie group $G$ defining $M$ is of type $A_l$ or $C_l$. There are many other cases, where equality holds, but in general not (see \cite{Diss}).

\vspace{0.2cm}

Concerning $l$: let $-K_M = rH$, $H$ ample. Then $H$ is very ample, since $M$ is rational homogeneous by \cite{RaRa}. Consider the embedding $\varphi_{|H|}: M \hookrightarrow \PN(H^0(M, H))$. Take exterior powers of the natural surjection $\Omega_{\PN}^1|_M \to \Omega_M^1 \to 0$ and use the duality $\bigwedge^lT_M \otimes K_M \simeq \Omega_M^{n -l}$. Then $\bigwedge^lT_M$ is a quotient of an ample bundle, hence ample, if 
 \begin{equation} \label{ineql}
   l \le n - r +2.
 \end{equation}
Let on the other hand $-K_M = \sum_{i=1}^{\rho(M)}k_iH_i$ in $\Pic(M)$, where $H_i$ are the fundamental divisors. Then $k_i \ge 1$ for all $i$. Let $k_{i_0}$ be a mimimal coefficient. Let $C$ be an $\alpha$-line in $M$ corresponding to that $i_0$, i.e., a smooth rational curve such that $C.H_j = \delta_{i_0j}$ (see \cite{Kollar}). Then $-K_M.C = k_{i_0}$. Let $(a_1, \dots, a_n)$ be the splitting type of $T_M$ on $C$. Since $T_M$ is globally generated, all $a_i \ge 0$. Moreover, $a_i \ge 2$ for at least one $i$, since $T_M|_C$ contains the tangent bundle of $C$. Then there are at least $n - \sum_i a_i +1 = n-k_{i_0}+1$ zeroes in the splitting type. That means, that $\bigwedge^jT_M$ cannot be ample for $j \le n-k_{i_0}+1$. Hence $l \ge n - \min_i(k_i) + 2$. Since $r = \gcd(k_i) \le \min_i(k_i)$, we get 
 \begin{equation} \label{minki}
   l = n - \min\nolimits_i(k_i) +2.
 \end{equation} 
If in particular $\rho(M) = 1$, then $\min_i(k_i) = r$, and we have equality in (\ref{ineql}).
\end{remark}


\section{Projective space and complex tori}
\setcounter{lemma}{0}

We will start with some comments on Van de Ven's theorem

\begin{theorem2}{Van de Ven}
 The tangent sequence of a submanifold of the projective space splits, if and only if it is a linear subspace. 
\end{theorem2}

Van de Ven himself published no proof of his theorem, but there are many different ones in the literature, see for example \cite{Laksov}, \cite{MoRo}, or \cite{Popa}. Working on projective connections, I. Radloff and the author found a very short proof in \cite{Prep}, using that $X$ inherits a projective connection from $\PN_n$. A direct proof can be given using Mori's theorem, that $\PN_n$ is the only projective manifold with ample tangent bundle (\cite{Mori}).

Another possibility is to use Wahl's characterization of $\PN_n$: if $X \subset \PN_n$ is a submanifold with splitting tangent sequence, then $T_X \otimes \sO_{\PN_n}(-1)|_X$ is globally generated, implying $(X, \sO_{\PN_n}(1)|_X) = (\PN_m, \sO_{\PN_m}(1))$ or $(\PN_1, \sO_{\PN_1}(2))$ by \cite{Wahl}. In the latter case use for example the Euler sequence to show that the tangent sequence of a conic in $\PN_2$ does not split. 

\vspace{0.2cm}

Next we consider compact complex manifolds with trivial tangent bundle, called parallelisable. A parallelisable K\"ahler manifold is a complex torus, in the projective case it is an abelian variety.

\begin{theorem} \label{Tori}
The tangent sequence of a submanifold of a parallelisable manifold splits, if and only if it is parallelisable.
\end{theorem}

\begin{proof} 
The proof relies on the following general fact on vector bundles (see for example \cite{OSS}, \S~5, Exercise): a globally generated vector bundle $E$ with trivial determinant on some manifold $X$ is trivial. Indeed: since $E$ is globally generated, we have a surjective map from a trivial bundle $\sO_X^{\oplus r}$ onto $E$. On the other hand, $E^* \simeq \bigwedge^{\rk(E)-1}E$ is globally generated as well, in particular has a section. This gives a map $E \to \sO_X$. We obtain a non-zero map $\sO_X^{\oplus r} \to \sO_X$, which is hence surjective and admits a splitting. We have proved $E \simeq E' \oplus \sO_X$. Moreover, $E'$ is again globally generated with trivial determinant. We proceed by induction on the rank. 

Let $X$ be a submanifold of the parallelisable manifold $M$ with splitting tangent sequence. Then $T_X$ is globally generated, so we just have to prove that $\det(T_X) = -K_X$ is trivial. Obviously, $-K_X$ is globally generated. On the other hand, $K_X$ is globally generated by adjunction formula, implying that $K_X$ is in fact trivial.

If conversely $X \subset M$ is parallelisable, then the normal bundle of $X$ in $M$ is globally generated with trivial determinant, hence trivial, i.e. the tangent sequence is a sequence of trivial bundles, hence splits.
\end{proof}


\section{Quadrics}
\setcounter{lemma}{0}

Let $Q_n \subset \PN_{n+1}$ be the $n$-dimensional quadric hyperplane. We will show that the list of submanifolds with splitting tangent sequence consists (almost) of what one expects from Van de Ven's theorem: linear subspaces and complete intersections of hyperplanes. Only curves match not exactly in that scheme, they are all rational, but may have higher degree (see Example~\ref{deg}). 

\vspace{0.2cm}

For $n \ge 3$, the quadric $Q_n$ is rational homogeneous with $\rho = 1$. The generator of $\Pic(Q_n)$ will be denoted by $\sO_{Q_n}(1) = \sO_{\PN_{n+1}}(1)|_{Q_n}$; by adjunction formula, $\sO_{Q_n}(-K_{Q_n}) = \sO_{Q_n}(n)$; by \cite{Gold}, $T_{Q_n}$ is $1$-ample. From the surjective map $\Omega_{\PN_{n+1}}^1 \to \Omega_{Q_n}^1 \to 0$ and the Euler sequence on $\PN_{n+1}$ we deduce
  \begin{equation} \label{Euler}    
    \sO_{Q_n}(1)^{\oplus N} \lra \bigwedge\nolimits^2 T_{Q_n} \lra 0
  \end{equation}
for $N = {n+2 \choose 3}$. The splitting type of $T_{Q_n}$ on a line $l$ in $Q_n$ is
  \begin{equation} \label{tblines}
    T_{Q_n}|_l = \sO_l \oplus \sO_l(1)^{\oplus (n-2)} \oplus \sO_l(2).
  \end{equation}
In particular, $\Seq_{l/Q_n}$ splits. The section in $\PN(T_{Q_n})$ of any line $l$ induced by the projection 
 \[T_{Q_n}|_l \lra \sO_l \lra 0\] 
will be contracted to a point by the morphism $\Phi_{Q_n}$ from (\ref{amp}). Using the identity $T_{Q_n} \simeq \Omega^1_{Q_n}(2)$, we may on the other hand view $\PN(T_{Q_n})$ as a subvariety of the flag variety of points and lines of $\PN_{n+1}$, such that the map $\Phi_{Q_n}$ is the projection to the Grassmannian of lines. Hence the exceptional fibers of $\Phi_{Q_n}$ are exactly the sections of lines in $Q_n$ (see \cite{Wis}, Example~1).

\vspace{0.2cm}

\noindent We start with $Q_2 \simeq \PN_1 \times \PN_1$:

\begin{proposition} \label{Q2}
 The tangent sequence of a submanifold of $Q_2 = \PN_1 \times \PN_1$ splits if and only if it is a rational curve.
\end{proposition}

\begin{proof}
 Let $C \subset Q_2$ be a smooth curve with splitting tangent sequence. Then $C$ is homogeneous, hence rational or elliptic. Assume $C$ elliptic. By adjunction formula, $C$ is a divisor of bidegree $(2,2)$ on $\PN_1 \times \PN_1$. Then the restriction of $T_{Q_2}$ to $C$ is the sum of two line bundles of degree $4$, whereas $T_C$ is trivial. A contradiction.

Let conversely $C$ be a rational curve in $Q_2$. Then $C$ is of type $(1, d-1)$ for some $d \ge 1$, i.e. $T_{Q_2}|_C = \sO_C(2) \oplus \sO_C(2d-2)$, proving the claim.
\end{proof}

\begin{proposition} \label{homratQ}
A submanifold with splitting tangent sequence in $Q_n$ is rational homogeneous.
\end{proposition}

\begin{proof}
We have already seen the case $n = 2$, hence assume $n \ge 3$. Let $X \subset Q_n$ be a submanifold with splitting tangent sequence. Then $X$ is homogeneous, hence $X \simeq A \times G/P$, $A$ abelian, $G/P$ rational homogeneous. Identify a fiber of the projection $X \to G/P$ with $A$. Then $\Seq_{A/Q_n}$ splits by \ref{trans}. The induced section $A'$ of $A$ in $\PN(T_{Q_n})$ will be contracted to a point by $\Phi_{Q_n}$, since $T_A$ is trivial. This is impossible if $\dim(A) > 0$, since all exceptional fibers of $\Phi_{Q_n}$ are rational curves. 
\end{proof}

\begin{abs} \label{curvesQ}
{\bf Rational curves of low degree in $Q_n$.} Let $C \subset Q_n$, $n\ge 3$, be a rational curve of degree $d \le 3$. Let $L \subset \PN_{n+1}$ be the smallest linear subspace containing $C$. Then $\dim(L) = d$ and $C$ is projectively normal embedded by the full linear system $\sO_C(d)$ (which is in general not true for $d \ge 4$), meaning the restriction map $H^0(L, \sO_L(1)) \to H^0(C, \sO_C(d))$ is an isomorphism. Therefore the dual Euler sequence
 \[0 \lra \Omega_L^1(1) \lra \sO_L^{\oplus (d+1)} \lra \sO_L(1) \lra 0\]
remains surjective on $H^0$-level on $C$. Hence $h^0(C, \Omega_L^1|_C(d)) = 0$. On the other hand, $T_L|_C(-d) = T_L(-1)|_C$ is globally generated. This shows $T_L|_C = \sO_C(d+1)^{\oplus d}$. Since $\Seq_{L/\PN_{n+1}}$ splits, we get $N_{C/\PN_{n+1}} = N_{C/L} \oplus \sO_C(d)^{\oplus (n+1-d)}$ and $T_{\PN_{n+1}}|_C = \sO_C(d+1)^{\oplus d} \oplus \sO_C(d)^{\oplus (n+1-d)}$.

Assume now the tangent sequence of $C$ in $Q_n$ splits. This is always the case if $C$ is a line (for $d \ge 2$ see the following Example~\ref{deg}). Let $(a_1, \dots, a_{n-1})$ be the splitting type of $N_{C/Q_n}$. Then $a_i \le d+1$ for all $i$ and $\oplus_{a_i = d+1}\sO_C(a_i) \hookrightarrow N_{C/L}$. Since the cokernel of this injective map is free and it cannot be an isomorphism since $\deg(N_{C/L}) = (d-1)(d+2)$, there are at most $\rk(N_{C/L})-1 = d-2$ of the $a_i$'s greater than $d$. Summing up the $a_i$'s, we find the splitting type of $T_{Q_n}|_C$ is
 \[(0, 1, \dots, 1, 2), \quad (2, \dots, 2) \quad \mbox{ and } \quad (2,3,\dots,3,4)\]
for $d = 1, 2$ and $3$, respectively. 
\end{abs}

\begin{example} \label{deg}
For $n \ge 6$ we construct rational curves of given degree $d \ge 2$ in $Q_n$, with splitting tangent sequence and with non-splitting one: consider two submanifolds with splitting tangent sequence in $Q_n$: a smooth $2$-dimensional complete intersection subquadric, and a $3$-dimensional linear subspace (which exists for $n \ge 6$). Both of them contain rational curves of arbitrary degree $d$. The tangent sequence of any rational curve in the subquadric splits, whereas the tangent sequence of a rational curve of degree $d \ge 2$ in the linear subspace does not.
\end{example}

\begin{theorem} \label{Qn}
The tangent sequence of a submanifold $X$ of $Q_n$ of dimension $m \ge 2$ splits if and only if $X \simeq \PN_m$ is a linear subspace or $X \simeq Q_m$ is a complete intersection subquadric. Curves in $Q_n$ with splitting tangent sequence are rational.
\end{theorem}

\begin{proof}
If $X \simeq\PN_m$ is a linear subspace of $Q_n$, then $\Seq_{X/\PN_{n+1}}$ splits by Van de Ven's theorem, hence $\Seq_{X/Q_n}$ splits by \ref{trans}. Let $Q_m \subset Q_n$ be a complete intersection subquadric, i.e. $N_{Q_m/Q_n} = \sO_{Q_m}(1)^{\oplus (n-m)}$. By \cite{SnowQ}, $H^1(Q_m, T_{Q_m}(-1)) = 0$, hence $\Seq_{Q_m/Q_n}$ splits by $h^1$-criterion.

Let conversely $X \subset Q_n$ be a submanifold with splitting tangent sequence. Then $X$ is rational homogeneous by Proposition~\ref{homratQ}. In particular, if $\dim(X) = 1$, then $X$ is a rational curve. Assume $\dim(X) = m \ge 2$. From the splitting $T_{Q_n}|_X \to T_X \to 0$ we obtain that $\bigwedge^2T_X$ is ample, using (\ref{Euler}). It follows $X \simeq \PN_m$ or $Q_m$: for $m = 2$ this is clear, since $\PN_2$ and $\PN_1 \times \PN_1$ are the only rational homogeneous surfaces. For $m \ge 3$, Cho-Sato's result applies (\cite{ChoSato}): the projective space and the quadric are the only projective manifolds $X$, such that $\bigwedge^2T_X$ is ample.

To complete the proof we have to determine the embeddings. Define $d = \deg(X)$ by $\sO_{Q_n}(1)|_X = \sO_X(d)$. Let $l \subset X$ be a line in $X$. Then $\Seq_{l/X}$ splits and $\sO_{Q_n}(1).l = d$. From \ref{trans} we get: $\Seq_{l/Q_n}$ splits, and the normal bundle sequence
 \[0 \lra N_{l/X} \lra N_{l/Q_n} \lra N_{X/Q_n}|_l \lra 0\]
splits as well. Let $(a_{m+1}, \dots, a_n)$ be the splitting type of $N_{X/Q_n}$ on $l$. The splitting type of $N_{l/X}$ is either $(1, \dots, 1)$, if $X \simeq \PN_m$, or $(0, 1, \dots, 1)$, in the case $X \simeq Q_m$. Denote the splitting type of $T_{Q_n}|_l$ by $(t_1, \dots, t_n)$. Then
 \[(t_1, \dots, t_n) = (1, \dots, 1, 2, a_{m+1}, \dots, a_n) \; \mbox{ or } \; (0, 1, \dots, 1, 2, a_{m+1}, \dots, a_n),\]
if $X \simeq \PN_m$ or $Q_m$, respectively. In \ref{curvesQ} we have proved: $(t_1, \dots, t_n) = (2, \dots, 2)$ and $(2, 3, \dots, 3, 4)$ for $d = 2,3$, respectively. This shows $d \not= 2,3$. The map from (\ref{Euler}) gives $\sO_l(d)^{\oplus N} \to \oplus_{i\not= j}\sO_l(t_i + t_j) \to 0$, hence $d \le t_i + t_j$ for all $i \not= j$. Since $m \ge 2$, we have $t_1 + t_2 \le 3$. The only remaining possibility is $d = 1$. If $X \simeq \PN_m$, we are already done, for $X \simeq Q_m$, we have to show that $X$ is a complete intersection of hyperplanes. 

Let $i: X = Q_m \hookrightarrow Q_n$ be the embedding, and consider $Q_m \subset \PN_{m+1}$ as a hypersurface of degree $2$. Since $H^0(Q_m, \sO_{Q_m}(1)) \simeq H^0(\PN_{m+1}, \sO_{\PN_{m+1}}(1))$, and 
 \[\sO_{Q_m}(1) = \sO_{Q_n}(1)|_{Q_m} = \sO_{\PN_{n+1}}(1)|_{Q_m},\] 
it is possible to extend $i$ to $\PN_{m+1}$. Hence $Q_m \subset \PN_{m+1} \subset \PN_{n+1}$, where the last $\PN_{n+1}$ is the one containing $Q_n$ as a hypersurface, and $\PN_{m+1} \subset \PN_{n+1}$ is a linear subspace. This linear subspace $\PN_{m+1}$ cannot be contained in $Q_n$: if it were, then the splitting of $\Seq_{Q_m/Q_n}$ would imply the splitting of $\Seq_{Q_m/\PN_{m+1}}$, which is not possible. It follows that $Q_m$ is the complete intersection of $Q_n$ and the linear subspace $\PN_{m+1}$ in $\PN_{n+1}$.
\end{proof}

\begin{remark} \label{alt}
Thinking of a generalization to other rational homogeneous manifolds, it might be interesting to have an alternative argument, avoiding Cho-Sato's beautiful characterization of $Q_m$. Since we already know, that $X$ is rational homogeneous by Proposition~\ref{homratQ}, we may use classification instead: let $-K_X = \sum_ik_iH_i$, $H_i \in \Pic(X)$ the fundamental divisors. Then the ampleness of $\bigwedge^2T_X$ implies $\min_i(k_i) \ge \dim(X)$ by (\ref{minki}) in Remark~\ref{kl}. Now check the list for the $k_i$, which can be found for example in \cite{Akh}, to verify that in fact $X = \PN_m, Q_m$ are the only rational homogeneous manifolds with this property.
\end{remark}


\section{Grassmannians and $\PN(T_{\PN_2})$}
\setcounter{lemma}{0}

The Grassmannian $\G = {\rm Gr}(k, n+1)$ of $k$-dimensional subspaces of $\KC^{n+1}$ is rational homogeneous with Picard number one. The dimension of $\G$ is $k(n+1-k)$, and by \cite{Gold}, $T_{\G}$ is $(k(n+1-k)-n)$-ample.

\begin{abs} \label{ampG}
{\bf The fibers of the morphism $\Phi_{\G}$.} There is a nice description of the fibers in this case: denote the trivial bundle on $\G$ by ${\mathcal V} = \sO_{\G}^{\oplus n+1}$. Let ${\mathcal U} \subset {\mathcal V}$ be the universal subsheaf and ${\mathcal Q}$ the universal quotient sheaf. Then $T_{\G} = \GHom({\mathcal U}, {\mathcal Q})$, and the projective bundle $\PN(T_{\G})$ consists of all pairs $(U, [\varphi])$ with $U\subset V$, $\dim(U) = k$, $\varphi\in\Hom(V/U,U)$ and $\varphi\not\equiv 0$. Thus, by composing inclusion and projection, there is a natural map 
 \[{\bf P}(T_{\G}^*) \stackrel{j}{\hookrightarrow} {\bf P}(\GHom(V,V)) \simeq \G \times {\bf P}(\Hom(V,V)) \stackrel{f}{\lra} {\bf P}(\Hom(V,V)).\]
Hence $\Phi_{\G} = f \circ j$ and the fiber over a point $[\varphi_0] \in {\bf P}(\Hom(V,V))$ is just the set of all pairs $(U, [\varphi]) \in {\bf P}(\GHom({\mathcal Q}, {\mathcal U}))$, such that $j(\varphi) = \varphi_0$, i.e., all $U \in \G$ with $\im(\varphi_0) \subset U \subset \ker(\varphi_0)$. The exceptional fibers of $\Phi_{\G}$ are isomorphic to ${\rm Gr}(k-r, n+1-2r)$, $1 \le r \le k-1$. This implies:
\end{abs}

\begin{proposition} \label{homratG}
 Let $\G = {\rm Gr}(k, n+1)$ be a Grassmann manifold. If $X \subset \G$ is a submanifold with splitting tangent sequence, then $X$ is rational homogeneous.
\end{proposition}

\begin{proof}
Let $X \subset \G$ be a submanifold with splitting tangent sequence. Then $X$ is homogeneous, hence $X \simeq  A \times G/P$, $A$ abelian and $G/P$ rational homogeneous. Assume $\dim(A) > 0$. Identify a fiber of the projection $X \to G/P$ with $A$. Then $\Seq_{A/\G}$ splits by \ref{fiberbundle} and \ref{trans}. The surjection 
 \[T_{\G} \to T_A = \sO_A^{\dim(A)}\] 
induces a section $A'$ of $A$ in $\PN(T_{\G})$, which is contracted to a point by $\Phi_{\G}$. But all exceptional fibers of $\Phi_{\G}$ are Grasmannians as seen in \ref{ampG}, meaning that there exists a strictly smaller Grassmannian $\G' \subset \G$ that contains $A$. By \ref{trans}, $\Seq_{A/\G'}$ splits. We proceed by induction on the dimension of $\G$. Since the smallest Grassmannian is a $\PN_1$, we finally get $\dim(A) = 0$.
\end{proof}

\

The flag manifold $\PN(T_{\PN_2})$ is rational homogeneous with Picard number two. It may be realized as a divisor of bidegree $(1,1)$ in $\PN_2 \times \PN_2$. The projections are denoted by $\pi_1, \pi_2: \PN(T_{\PN_2}) \to \PN_2$; the two fundamental divisors are 
 \[L_i = \pi_i^*\sO_{\PN_2}(1), \quad i = 1,2.\] 
We consider here just the case $n = 2$ instead of $\PN(T_{\PN_n})$ in general, since we then may determine again the fibers of $\Phi_{\PN(T_{\PN_2})}$: similar to the case of $Q_n$ we find that any exceptional fiber of $\Phi_{\PN(T_{\PN_2})}$ is a section of some $F = F_1 + F_2$ in $\PN(T_{\PN_2})$, where $F_i$ is a fiber of the projection $\pi_i$ (see \cite{Wis}, Example~2). Now the same argument as in Proposition~\ref{homratG} shows:

\begin{proposition} \label{homratPvT}
If $X \subset \PN(T_{\PN_2})$ is a submanifold with splitting tangent sequence, then $X$ is rational homogeneous.
\end{proposition}

\begin{proposition} \label{PvT} 
The tangent sequence of a submanifold of $\PN(T_{\PN_2})$ splits if and only if it is either a section of a line in $\PN_2$ or a rational curve $C$ of bidegree $(2,2)$, such that $\pi_i|_C: C \to \PN_2$ is an embedding for both projections.
\end{proposition}

\begin{proof}
1.) The tangent sequence of a section of a line in $\PN_2$ splits by \ref{projbundle}. Let $C \subset \PN(T_{\PN_2})$ be a rational curve of bidegree $(2,2)$, such that $\pi_i|_C: C \to \PN_2$ is an embedding for $i = 1,2$. 

First note that such curves really exist: take any smooth conic $C_0$ in $\PN_2$. Since the splitting type of $T_{\PN_2}$ on $C_0$ is $(3,3)$, we have $S_1 = \pi_1^{-1}(C_0) \simeq \PN_1 \times \PN_1$, the second projection $\pi_2: S_1 \to \PN_2$ being a double cover. The ramification divisor is our curve $C \in |\sO_{\PN_1 \times \PN_1}(1,1)|$. Conversely, $C$ will always be contained in the surfaces 
 \[S_i = \pi_i^{-1}(\pi_i(C)) \simeq \PN_1 \times \PN_1,\] 
since $\pi_i(C) \subset \PN_2$ are smooth conics by assumption. 

We want to show the splitting of $\Seq_{C/\PN(T_{\PN_2})}$. The restriction of the relative tangent sequence corresponding to $\pi_i$ is
 \[0 \lra \sO_{\PN_1}(2) \lra T_{\PN(T_{\PN_2})}|_C \lra \sO_{\PN_1}(3) \oplus \sO_{\PN_1}(3) \lra 0.\]
This sequence splits, i.e., the splitting type of $T_{\PN(T_{\PN_2})}|_C$ is $(2,3,3)$. We claim the splitting type of the normal bundle is $(3,3)$. Consider
 \[\alpha: X \lra \PN(T_{\PN_2})\]
the blow-up of $\PN(T_{\PN_2})$ along $C$. By classification, $X$ is a Fano threefold with $\rho(X) = 3$ (see \cite{MoriMukaiM} and \cite{MoriMukaiA} or \cite{AG5}), and we have the following commutative diagram (\cite{MoriMukaiM}, Proposition~9 and 10)
 
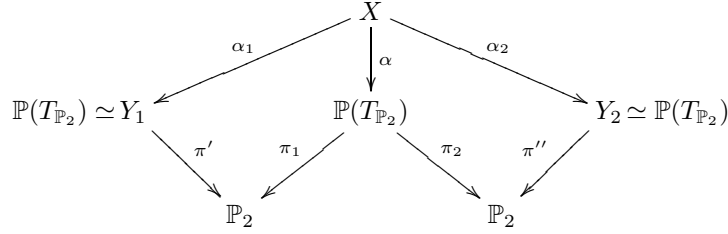
\[\xymatrix{& & & X \ar[lld]_{\alpha_1} \ar[rrd]^{\alpha_2} \ar[d]^{\alpha} & & &\\
  \PN(T_{\PN_2}) \simeq \hspace{-1cm} & Y_1 \ar[dr]^{\pi'} & & \PN(T_{\PN_2}) \ar[dl]_{\pi_1} \ar[dr]^{\pi_2} & & Y_2 \ar[dl]_{\pi''} & \hspace{-1cm} \simeq \PN(T_{\PN_2})\\         
    & & \PN_2 & & \PN_2 & &}\]
where $\alpha_i$ is birational for $i = 1$ and $2$, contracting the strict transform of the surface $S_i$ to a curve $C_i$ of the same type as $C$, i.e., $C_i$ is a curve of bidegree $(2,2)$, such that both projections embed $C_i$ into $\PN_2$. Starting now with say $C_1$ in $Y_1 \simeq \PN(T_{\PN_2})$, we find that the exceptional divisor 
 \[E = \PN(N^*_{C/\PN(T_{\PN_2})})\] 
of the blow-up $\alpha$ is the strict transform of $S' = {\pi'}^{-1}(\pi'(C_1)) \simeq \PN_1 \times \PN_1$. This shows $E \simeq \PN_1 \times \PN_1$. By adjunction formula, $\deg(N_{C/\PN(T_{\PN_2})}) = 6$, i.e., the splitting type of $N_{C/\PN(T_{\PN_2})}$ is indeed $(3,3)$, proving the splitting of $\Seq_{C/\PN(T_{\PN_2})}$.

\vspace{0.2cm}

2.) For the converse we have to consider rational curves, and surfaces isomorphic to $\PN_2$ or $\PN_1 \times \PN_1$ by Proposition~\ref{homratPvT}.
 
(i) Let $C \subset \PN(T_{\PN_2})$ be a rational curve with splitting tangent sequence. Assume $L_1.C = a$ and $L_2.C = b$ with $a \ge b \ge 0$. If $b = 0$, then $C$ is a fiber of $\pi_2$ and we are done. If $b = 1$, then $\pi_2$ maps $C$ isomorphically to a line in $\PN_2$. We may hence assume $b \ge 2$.

Assume $b = 2$. Then $\pi_2^*\sO_{\PN_2}(1).C = 2$ implies $\pi_2|_C$ is either a double cover of a line or an isomorphism onto a smooth conic in $\PN_2$. Assume first $\pi_2(C)$ is a line. Then $C$ is contained in the Hirzebruch surface 
 \[\Sigma = \PN(T_{\PN_2}|_{\pi_2(C)}) = \PN(\sO_{\PN_1}(1) \oplus \sO_{\PN_1}(2))\] 
with $C \in |\sO_{\Sigma}(2) \otimes \pi_2|_{\Sigma}^*\sO_{\PN_1}(-2)|$, since $C$ is rational by assumption. We find that $C$ does not meet the minimal section of $\Sigma$, i.e., viewing $\Sigma = \PN_2(p)$ as blow-up of $\PN_2$ in a point, $\Sigma$ is isomorphic to that $\PN_2$ in a neighborhood of $C$. But the splitting of $\Seq_{C/\PN(T_{\PN_2})}$ implies the splitting of $\Seq_{C/\Sigma}$, which is impossible, since the tangent sequence of a conic in $\PN_2$ does not split.

Assume $b = 2$ and $\pi_2(C)$ is a conic in $\PN_2$. This is impossible for $a \ge 3$, since the tangent sequence of a conic in $\PN_2$ does not split, and \ref{projbundle} applies. Hence $a = 2$ and as above we find $\pi_1|_C: C \to \PN_2$ is an embedding as well.

Assume $b \ge 3$. Let $(s, t)$ be the splitting type of $\pi_2^*T_{\PN_2}|_C$. Then $s,t \ge b$ (Euler sequence). For the relative tangent sequence corresponding to $\pi_2$ we see $T_{\PN(T_{\PN_2})/\PN_2}|_C = \sO_C(2a-b)$. The dual sequence twisted by $\sO_C(2)$ then yields $2-s \ge 0$ or $2-t \ge 0$, since $b-2a+2 < 0$ (use $a \ge b \ge 3$), a contradiction.

\vspace{0.2cm}

(ii) Let $S \subset \PN(T_{\PN_2})$ be a surface with splitting tangent sequence. From the ideal sequence we obtain $h^0(\PN(T_{\PN_2}), \sO_{\PN(T_{\PN_2})}(K_{\PN(T_{\PN_2})} + S)) = 0$ since 
 \[h^0(\PN(T_{\PN_2}), K_{\PN(T_{\PN_2})}) = h^0(S, K_S) = 0.\]
Assume $S \in |aL_1+bL_2|$ with $a \ge b \ge 0$. Consider $\PN(T_{\PN_2})$ as a divisor of bidegree $(1,1)$ in $\PN_2 \times \PN_2$. Then the above vanishing implies 
 \[h^0(\PN_2 \times \PN_2, \sO_{\PN_2 \times \PN_2}(a-2,b-2)) = h^0(\PN_2 \times \PN_2, \sO_{\PN_2 \times \PN_2}(a-3,b-3)),\]
which implies $b \le 1$ by K\"unneth's formula. If $b = 0$, then there exists a smooth curve $C$ of degree $a$ in $\PN_2$, such that $S \simeq \PN(T_{\PN_2}|_C)$. Then $C$ is rational and $S \simeq \PN_1 \times \PN_1$, which is impossible. If $b = 1$, then $S \simeq \PN_2$ is a section of $\pi_1$ and $9 = K_S^2 = 8-a(a+1)$ by adjunction formula, again a contradiction.
\end{proof}


\section{Blow-up of $\PN_n$ in a point}
\setcounter{lemma}{0}

The blow-up of $\PN_n$ in a point $p \in \PN_n$ is not homogeneous but almost homogeneous. We use the classification theorem of Huckleberry and Oeljeklaus for almost homogeneous varieties with exceptional locus a point (\cite{HuOe}). 

\vspace{0.2cm}

Denote the blow-up of $\PN_n$ in the point $p$ by $\varphi: \PN_n(p) \to \PN_n$. Call the exceptional divisor $E$. Then $|\varphi^*\sO_{\PN_n}(1)-E|$ defines a morphism $\pi: \PN_n(p) \to \PN_{n-1}$, identifying $\PN_n(p)$ with the projective bundle $\PN(\sO_{\PN_{n-1}}\oplus\sO_{\PN_{n-1}}(1))$:
 \[\xymatrix{\PN_n(p) \ar[r]^{\pi} \ar[d]_{\varphi} & \PN_{n-1}\\
   \PN_n &}\]
The tautological line bundle $\sO(1)$ of $\PN(\sO_{\PN_{n-1}}\oplus\sO_{\PN_{n-1}}(1))$ is $\varphi^*\sO_{\PN_n}(1)$, the exceptional divisor $E$ of $\varphi$ is the minimal section of $\pi$. We have

\begin{lemma} \label{point}
Let $\PN_m \subset \PN_n$ be a linear subspace. Then the tangent sequence of the strict transform of $\PN_m$ in $\PN_n(p)$ splits. 
\end{lemma}

\begin{proof}
If $\PN_m$ does not contain the point $p$, then the assertion is trivial, we may hence assume $p \in \PN_m$. Then the strict transform of $\PN_m$ is $\PN_m(p)$. By the transitivity principle \ref{trans}, we may assume $n = m+1 \ge 2$. The claim hence follows from the vanishing of
 \[H^1(\PN_m(p), T_{\PN_m(p)} \otimes \pi^*\sO_{\PN_{m-1}}(-1)).\] 
We consider the relative tangent sequence
 \[0 \lra \sO(2) \otimes \pi^*\sO_{\PN_{m-1}}(-1) \lra T_{\PN_m(p)} \lra \pi^*T_{\PN_{m-1}} \lra 0,\]
twisted by $\pi^*\sO_{\PN_{m-1}}(-1)$. Since all higher direct image sheaves vanish, we obtain the desired vanishing for $m \ge 3$ by Bott's formula. For $m=2$, consider the twisted tangent sequence
 \[0 \lra T_{\PN_2(p)} \otimes \pi^*\sO_{\PN_1}(-1) \lra T_{\PN_3(p)} \otimes \pi^*\sO_{\PN_2}(-1)|_{\PN_2(p)} \lra \sO_{\PN_2(p)} \lra 0.\]
We have to show that this sequence is exact on $H^0$-level, which follows from the identity 
 \[H^1(\PN_2(p), T_{\PN_2(p)} \otimes \pi^*\sO_{\PN_1}(-1)) \simeq H^1(\PN_2(p), T_{\PN_3(p)} \otimes \pi^*\sO_{\PN_1}(-1)|_{\PN_2(p)}),\] 
since $H^1(\PN_2(p), \sO_{\PN_2(p)})$ vanishes. Applying the above argument twice, we show that both $H^1$'s are onedimensional, where we use the twisted ideal sequence of $\PN_2(p)$ in $\PN_3(p)$ to compute the right hand side. In the case $m=1$, $\PN_1(p) \simeq \PN_1$ is a fiber of $\pi$, i.e. the claim follows by \ref{fiberbundle}.
\end{proof}

\

Now let $X \subset \PN_n(p)$ be a submanifold with splitting tangent sequence, not contained in $E$. There are vector fields on $\PN_n$ with zero in $p$, generating $T_{\PN_n}$ outside $p$. These generate $T_{\PN_n(p)}$ outside $E$ and the same holds true for the quotient bundle $T_X$. Therefore $X$ and $\varphi(X)$ are almost homogeneous varieties with exceptional locus contained in $E$ and $\{p\}$, respectively. We will start considering curves:

\begin{abs} \label{curvesPn(p)}
{\bf Curves in $\PN_n(p)$ with splitting tangent sequence.} Let $C \subset \PN_n(p)$ be a smooth curve. If $C$ is contained in $E$, then $\Seq_{C/\PN_n(p)}$ splits, if and only if $C$ is a line in $E$; if $C$ does not meet $E$, then $\Seq_{C/\PN_n(p)}$ splits, if and only if $C$ is the strict transform of a line in $\PN_n$. 

Assume $\Seq_{C/\PN_n(p)}$ splits and $C$ meets $E$ in points. Then there are vector fields with zeroes on $C$, hence $1 \le E.C \le 2$ and $C$ is rational. Define 
 \[\quad k = E.C \quad \mbox{ and } \quad d = \varphi^*\sO_{\PN_n}(1).C.\] 
We have $1 \le k \le 2$ and $d \ge k$. From $\pi^*\sO_{\PN_{n-1}}(1) = \varphi^*\sO_{\PN_n}(1) \otimes \sO_{\PN_n(p)}(-E)$ we obtain 
 \[\pi^*\sO_{\PN_{n-1}}(1).C = d-k.\] 
Show $d - k \le 2$. Let $(t_1, \dots, t_{n-1})$ be the splitting type of $\pi^*T_{\PN_{n-1}}|_C$. Pulling back the Euler sequence shows $t_i \ge d-k \ge 0$ for all $i$. Consider the relative tangent sequence restricted to $C$
 \begin{equation} \label{reltang} 
  0 \lra \sO_C(d+k) \lra T_{\PN_n(p)}|_C \lra \pi^*T_{\PN_{n-1}}|_C \lra 0.
 \end{equation}
From the splitting we have a ``2'' in the splitting type of $T_{\PN_n(p)}|_C$. Then either $d = k = 1$, or there is at least one $t_i \le 2$. Hence $0 \le d - k \le 2$ as claimed. We consider the cases $d-k = 0,1,2$ seperately, having in mind that $k = 1,2$.

\vspace{0.2cm}

1.) $d-k = 0$. Then $d = k = 1$, i.e. $C$ is a fiber of $\pi$, that is the strict transform of a line in $\PN_n$ and $\Seq_{C/\PN_n(p)}$ indeed splits by Lemma~\ref{point}.

\vspace{0.2cm}

2.) $d-k = 1$. We claim that $\Seq_{C/\PN_n(p)}$ splits if and only if we are in the situation of 3. in Theorem~\ref{Pn(p)}. The morphism $\pi|_C$ is defined by sections from $|\pi^*\sO_{\PN_{n-1}}(1)|_C| = |\sO_C(1)|$, hence $C$ maps isomorphically to a line in $\PN_{n-1}$. Then $C$ is a section of the Hirzebruch surface 
 \begin{equation} \label{sigma}
   \Sigma = \PN(\sO_{\PN_{n-1}}|_{\pi(C)} \oplus \sO_{\PN_{n-1}}(1)|_{\pi(C)}) = \PN(\sO_{\PN_1} \oplus \sO_{\PN_1}(1)).
 \end{equation}
The relative tangent sequence restricted to $C$ looks like
 \[0 \lra \sO_C(d+k) \lra T_{\Sigma}|_C \lra \sO_C(2) \lra 0\] 
and splits by $h^1$-criterion, i.e. $T_{\Sigma}|_C = T_C \oplus \sO_C(d+k)$. Then $\Seq_{C/\Sigma}$ splits. On the other hand, $\Sigma$ is the blow-up of a linear subspace $\PN_2 \subset \PN_n$ in $p$, meaning $\Seq_{\Sigma/\PN_n(p)}$ splits by Lemma~\ref{point}. This implies the splitting of $\Seq_{C/\PN_n(p)}$ (transitivity principle \ref{trans}).
The image $\varphi(C)$ is a plane curve of degree $d$ in $\PN_n$. There are two possibilities: if $k = 1$, then $d = 2$ and $\varphi(C)$ is a smooth conic in $\PN_n$; if $k = 2$, then $d = 3$ and $\varphi(C)$ is a curve of degree $3$ with a node in $p$. 

\vspace{0.2cm}

3.) $d-k = 2$. We will show that $\Seq_{C/\PN_n(p)}$ does not split in this case. Since $\pi^*\sO_{\PN_{n-1}}(1).C = 2$, $\pi|_C$ is either an isomorphism onto a smooth conic in $\PN_{n-1}$, or a 2:1-covering of a line. Assume the latter case. Then $C$ is a multisection of degree $2$ in the Hirzebruch surface $\Sigma$ from (\ref{sigma}). By adjunction formula, $C \in |\sO_{\Sigma}(2)|$, where $\sO_{\Sigma}(1)$ is the tautological line bundle on $\Sigma$. If $C_0$ denotes the minimal section of $\Sigma$, then $C.C_0 = 0$. On the other hand, $C_0$ is the restriction of $E$ to $\Sigma$, hence $C.C_0 = k$ by definition. This is impossible.

Assume $\pi|_C: C \to \pi(C)$ is an isomorphism. We have $T_{\PN_n(p)/\PN_{n-1}}|_C = \sO_C(d+k)$ from (\ref{reltang}) and $d + k = 2k+2 \ge 4$, i.e. $\Hom(T_{\PN_n(p)/\PN_{n-1}}|_C, T_C) = 0$. Since $\Seq_{\pi(C)/\PN_{n-1}}$ does not split, $\Seq_{C/\PN_n(p)}$ cannot split by \ref{projbundle}.
\end{abs}

\begin{theorem} \label{Pn(p)}
 The tangent sequence of a submanifold $X$ of the blow-up $\PN_n(p)$ of $\PN_n$ in a point $p$ splits if and only if one of the following holds:
 \begin {enumerate}
   \item $X$ is a linear subspace of the exceptional divisor,
   \item $X$ is the strict transform of a linear subspace in $\PN_n$,
   \item $X$ is the strict transform of a smooth conic through $p$ or of a plane curve of degree $3$ with node in $p$.
 \end{enumerate}
\end{theorem}

\begin{proof}
We have seen the result for curves in \ref{curvesPn(p)}, so assume $\dim X \ge 2$. Then Lemma~\ref{point} shows the ``if''-part of the theorem. It remains to prove the converse.

Let $X \subset \PN_n(p)$ be a submanifold with splitting tangent sequence of dimension at least $2$. If $X \cap E = \emptyset$ or $X \subset E$, then the claim is trivial by Van de Ven's theorem and \ref{trans}. We may hence assume $X \cap E \not\in \{\emptyset, X\}$. 

From the splitting, the image $\varphi(X) \subset \PN_n$ is an almost homogeneous variety of dimension greater than $1$ with exceptional locus contained in $\{p\}$. Then by \cite{HuOe}, p.113, {\em Theorem~1}, the normalization $\mu : X' \to \varphi(X)$ is a cone over a rational homogeneous manifold $Q$. That means: $X'$ is the image of the contraction of the minimal section $Q_0$ of the projective bundle $Z = \PN(\sO_Q \oplus L)$ over $Q$ for some very ample line bundle $L$ on $Q$:
 \[\xymatrix{
    && X \ar[dl]^{\varphi'} \ar[d]^{\varphi|_X} \ar[r]^{\pi|_X} & \PN_{n-1}\\  
    Z \ar[r]_{\psi} \ar[d]_{\rho} \ar@{-->}[urr]^{\chi} &  X' \ar[r]_{\mu} &  \varphi(X) &\\
    Q &&&&}\]
By construction, $X'$ is smooth outside $\psi(Q_0) = p'$ and $p'$ lies over $p$. The induced map $\varphi'$ is an isomorphism on $X\backslash (E\cap X) \stackrel{\sim}{\lra} X'\backslash \mu^{-1}(p)$, where $\mu^{-1}(p)$ is a finite set containing $p'$. On the other hand, $\psi$ is an isomorphism on $Z\backslash Q_0 \stackrel{\sim}{\lra} X'\backslash \{p'\}$, hence 
 \begin{equation} \label{points}
   \chi = {\varphi'}^{-1} \circ \psi\colon \;  Z\backslash (Q_0 \cup \mathrm{points}) \stackrel{\sim}{\lra} X\backslash (E\cap X)
 \end{equation}
is an isomorphism. Choose a general $Q_1 \in |\sO_{Z}(1)|$, the tautological line bundle on $Z$, such that $Q_1$ does not meet the exceptional points in (\ref{points}). Since $Q_1 \cap Q_0 = \emptyset$, we find an isomorphic image $\chi(Q_1) \simeq Q_1$ in $X$, not meeting $E$. Since $\Seq_{Q_1/Z}$ splits by \ref{projbundle}, also $\Seq_{\chi(Q_1)/X}$ splits. This implies the splitting of $\Seq_{\chi(Q_1)/\PN_n(p)}$ by \ref{trans}. Since $\chi(Q_1)$ does not meet $E$, the image $\varphi(\chi(Q_1))$ is isomorphic to $\chi(Q_1)$ in a neighborhood of $\chi(Q_1)$, hence $\varphi(\chi(Q_1))$ is smooth and $\Seq_{\varphi(\chi(Q_1))/\PN_n}$ again splits. By Van de Ven's theorem, then $\varphi(\chi(Q_1)) \simeq \PN_m$ is a linear subspace. 

We conclude: $Q \simeq \PN_m$ and $L = \sO_{\PN_m}(k)$ for some $k \ge 1$. Since $\psi$ is defined by $|\sO_Z(1)|$, we find $\psi^*\mu^*\sO_{\PN_n}(1) = \sO_Z(d)$ for some $d \ge 1$. Restricting to $Q_1 \simeq \PN_m$ gives $kd = 1$, since the image of $Q_1$ in $\PN_n$ is a linear subspace. Hence $k = d = 1$, i.e., $X' \simeq \PN_{m+1}$ is smooth and $\mu$ is an isomorphism. Then $\varphi(X) \simeq X' \simeq \PN_{m+1}$ is a linear subspace of $\PN_n$, completing the proof.  
\end{proof}


\end{document}